\documentclass[a4paper, 12pt, review]{elsarticle}
\usepackage[utf8]{inputenc}
\usepackage[margin=1in]{geometry}
\usepackage{amsmath}
\usepackage{amsfonts}
\usepackage{amssymb}
\usepackage{caption}
\usepackage{subcaption}
\usepackage{algorithm2e}
\usepackage{tikz}
    \usetikzlibrary[shapes.geometric] 
\usepackage[inline]{enumitem}
\usepackage[implicit=false]{hyperref}
\usepackage[theme=grayscale]{jlcode}

\makeatletter
\def\ps@pprintTitle{%
 \let\@oddhead\@empty
 \let\@evenhead\@empty
 \def\@oddfoot{}%
 \let\@evenfoot\@oddfoot}
\makeatother

\newcommand{\braces}[1]{\left\{#1\right\}}

\definecolor{blue}{rgb}{0.23,0.58,0.89}

\title{Solving influence diagrams via efficient mixed-integer programming formulations and heuristics}

\author{Helmi Hankimaa}
\author{Olli Herrala}
\author{Fabricio Oliveira\corref{cor1}}
\author{Jaan Tollander de Balsch}

\cortext[cor1]{Corresponding author: fabricio.oliveira@aalto.fi}

\address{Department of Mathematics and Systems Analysis,
Aalto University,
Espoo, Finland}

\begin{document}

\begin{abstract}
     In this paper, we propose novel mixed-integer linear programming (MIP) formulations to model decision problems posed as influence diagrams. We also present a novel heuristic that can be employed to warm start the MIP solver, as well as provide heuristic solutions to more computationally challenging problems. 
    We provide computational results showcasing the superior performance of these improved formulations as well as the performance of the proposed heuristic. Lastly, we describe a novel case study showcasing decision programming as an alternative framework for modelling multi-stage stochastic dynamic programming problems. 
\end{abstract}


\maketitle

\section{Introduction}

A powerful way to tackle complex real-life decision-making problems is to frame them as multi-stage decision problems under uncertainty. This helps one to understand the interactions between parts of the decision process and how uncertainty behaves from a rigorous, quantitative standpoint. State-of-the-art approaches for modeling and solving multi-stage decision problems stem from two main areas: decision analysis and stochastic programming. Each field has provided time-tested modeling and solution approaches with varying degrees of success depending on how uncertainty is modeled. 

These modeling approaches require specifying probability values and measurable consequences to uncertain events such that the decisions maximizing an expected utility function value can be made. However, it remains extremely challenging to analytically model decision problems due to the particularly interdependent nature that decision problems can exhibit. Outcomes of random events change the optimal decisions and decisions can influence the probability distributions of random events. The task of representing such dependencies may be equally challenging, if not more so, than finding the best decisions.

Influence diagrams \citep{howard2005influence} provide both a formal description of a decision problem and serve as a communication tool which requires minimal technical proficiency. Furthermore, they are useful in conveying structural relationships of the problem in a simple manner and thus crucially bridge the gap between quantitative specifications and qualitative descriptions. Due to their generality, influence diagrams pervade many modeling-based approaches that require a formal description of relationships between uncertainty, decisions and consequences. 

\begin{figure}[h]
    \centering
    \begin{tikzpicture}
    [decision/.style={fill=blue!80, draw, minimum size=2em, inner sep=2pt}, 
    chance/.style={circle, fill=orange!80, draw, minimum size=2em, inner sep=2pt},
    value/.style={diamond, fill=teal!80, draw, minimum size=2em, inner sep=2pt},
    scale=1.2]
    \node[chance]   (L) at (-0.25, 1.5)  {$L$};
    \node[chance]   (L1) at (0.75, 3)  {$R_1$};
    \node[chance]   (L2) at (0.75, 2)  {$R_2$};
    \node[chance]   (LN-1) at (0.75, 1)  {$R_{\dots}$};
    \node[chance]   (LN) at (0.75, 0)  {$R_N$};
    \node[decision] (1) at (1.75, 3)  {$A_1$};
    \node[decision] (2) at (1.75, 2)  {$A_2$};
    \node[decision] (N-1) at (1.75, 1)  {$A_{\dots}$};
    \node[decision] (N) at (1.75, 0)  {$A_N$};
    \node[chance]   (F) at (2.75, 1.5)  {$F$};
    \node[value]    (T) at (3.75, 1.5)  {$T$};     
    \draw[->, thick] (L) -- (L1);
    \draw[->, thick] (L) -- (L2);
    \draw[->, thick] (L) -- (LN-1);
    \draw[->, thick] (L) -- (LN);
    \draw[->, thick] (L) -- (F);
    \draw[->, thick] (L1) -- (1);
    \draw[->, thick] (L2) -- (2);
    \draw[->, thick] (LN-1) -- (N-1);
    \draw[->, thick] (LN) -- (N);
    \draw[->, thick] (1) -- (F);
    \draw[->, thick] (2) -- (F);
    \draw[->, thick] (N-1) -- (F);
    \draw[->, thick] (N) -- (F);
    \draw[->, thick] (F) -- (T);
    \draw[->, thick] (1) -- (T);
    \draw[->, thick] (N) -- (T);
\end{tikzpicture}
    \caption{An influence diagram representing the N-monitoring problem \cite{salo2022}. Decisions are represented by squares, chance events by circles and consequences by diamonds}
    \label{fig:n_monitoring}
\end{figure}
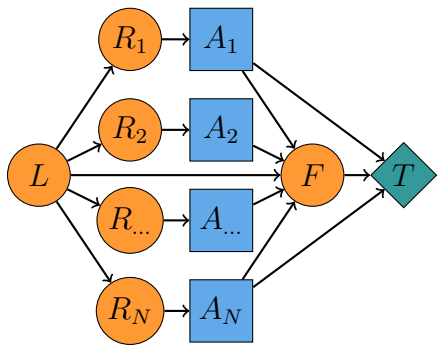

Figure~\ref{fig:n_monitoring} shows an influence diagram for the $N$-monitoring problem, where $N$ agents ($A_1$ to $A_N$) must decide whether to countervail an unknown load ($L$) based on imprecise readings of this load from their respective sensors ($R_1$ to $R_N$). The chance of failure ($F$) is influenced by the unknown load and the eventual decision to countervail the load. The final utility ($T$) is calculated considering whether the agents intervened and if a failure was observed. As such, this setting represents independent agents who must make decisions based only on observations of the state of the world but without being able to completely know it or share information among themselves.

Differently from decision trees, the nodes in an influence diagram do not need to be totally ordered nor do they have to depend directly on all predecessors. This freedom from dependence on all predecessors allows for the decisions to be made by, e.g., decision-makers who partially observe a common state of information (node $L$ in Figure~\ref{fig:n_monitoring}) but may differ in their ability to observe or are incapable of sharing information. 

Unfortunately, quantitative methods employed to obtain optimal decisions from influence diagram representations typically require that some of that generality is curbed. Influence diagrams are, in essence, representations of (possibly partially observable) Markov decision processes \citep{puterman1990markov}. Thus, if 
\begin{enumerate*}[label=(\roman*)]
    \item a single decision-maker is assumed (implying a total ordering among decision nodes) and\label{assumption:single_DM}
    \item the \emph{no-forgetting} assumption holds (implying that each decision node and its direct predecessors influence all successor decision nodes),\label{assumption:no_forgetting}
\end{enumerate*}
then Markovian assumptions hold. This, in turn, enables one to solve influence diagrams with well-established techniques, for instance, by forming the equivalent decision tree that can be solved through dynamic programming; or by removing decision and chance nodes from the diagram one by one \citep{bielza2011review}. The different issues regarding, e.g., computational performance of solving influence diagrams are discussed in \citet{bielza2000structural}.

As one may suspect, many problems, including that illustrated in Figure~\ref{fig:n_monitoring}, violate assumptions \ref{assumption:single_DM} and \ref{assumption:no_forgetting}. Indeed, there may be no memory or communication between deciding agents (meaning that they cannot know each other’s decisions) or constraints imposed across the diagram, such as budget limitations or logical conditions (e.g., stating that a given action can only occur if a prerequisite project has been started/completed in the past). All of these either violate the assumption that the previous state is ``remembered'' at a later stage or that all information influencing the decision alternatives is known when making said decision. These limitations create severe deficiencies in representing real-world problems. 

\citet{lauritzen2001representing} proposed an analytical framework to characterise these \emph{limited memory influence diagrams}. Note that the notion of limited memory can also be used to encompass settings with multiple decision-makers, the limited or absent sharing of information being its defining feature (regardless of whether it is due to lack of memory or communication between decision-makers). In any case, limited-memory influence diagrams are essentially influence diagrams that do not satisfy assumptions \ref{assumption:single_DM} and \ref{assumption:no_forgetting}. However, the fact that they do not satisfy assumptions \ref{assumption:single_DM} and \ref{assumption:no_forgetting} means that much more sophisticated analysis is required for designing methods that can extract optimal decisions from these diagrams. These methods involve specialised structures such as junction trees and message-passing mechanisms. 
  
Existing approaches for solving influence diagrams are deficient in many ways. First, solution methods are typically based on ad-hoc algorithm implementations that require the user to be proficient in manipulating influence diagrams. Additionally, such methods are often designed only for problems where expected utility is maximized, with no constraints connecting different decisions (e.g., budget/chance constraints). This limits their capabilities in modeling risk-averse decision making. Finally, tackling the challenges related to implementing computationally efficient methods (e.g., memory allocation, and thread-safe parallelisation) is required for larger problems. This creates a significant entry barrier to the wider adoption of such approaches.

Aiming to address the aforementioned limitations, Decision programming \citep{salo2022} leverages the capabilities of stochastic programming \citep{birge2011introduction} and decision analysis to model and solve multi-stage decision problems using mathematical optimization techniques. In essence, decision programming exploits the expressiveness of influence diagrams in structuring problems to develop deterministic equivalent \citep{birge2011introduction} mixed-integer programming (MIP) formulations. In turn, commercial-grade, professionally developed, off-the-shelf software can be readily utilised for solving these problems in a relatively straightforward manner instead of relying on ad-hoc implementations. It is worth highlighting that the latter, in stark contrast with the former, tend to be problem-specific and provide few guarantees for complying with sound software engineering practices in terms of versioning, updating and continuous improvements.

Concomitantly, the fact that decision programming models can be formulated as MIP problems gives rise to major benefits from the modeling perspective. In particular, decision programming exploits the exceptional modeling expressiveness provided by MIP to tackle challenging decision problems originally posed in the field of decision analysis. Having an underpinning MIP-based approach is extremely timely, especially due to the recent remarkable progress in MIP technology, with some estimating that the combined hardware and software speed-ups amount to a factor of two trillion between the early 1990s and mid-2010s \citep{bertsimas2019machine}! Thus, there is convincing evidence that MIP approaches for decision analysis problems are practically relevant, perhaps contrary to a few decades ago. 

In this paper, we provide multiple contributions in terms of practical and methodological aspects associated with decision programming. First and foremost, we present novel contributions that further improve upon the original ideas in \citet{salo2022}. Specifically, we present a novel formulation for decision programming problems that is considerably more efficient from a computational standpoint. This is illustrated in the computational experiments we provide. Furthermore, we propose a heuristic inspired by the single policy update heuristic, originally proposed in \citet{lauritzen2001representing}. This heuristic can be used not only to generate feasible solutions in case of more computationally challenging problems but also can be employed to warm start the MIP model.

Finally, we present a comprehensive case study in which we use the decision programming framework to develop an optimal decision strategy for allocating preventive care for coronary heart disease (CHD). The aim of this study is to evaluate the suitability of decision programming for performing the cost-benefit analysis originally performed by \citet{hynninen2019value}. In \citet{hynninen2019value}, a set of alternative predefined testing and treatment strategies for CHD are optimized using dynamic programming. We show how the same problem can be solved precluding the need to define strategies a priori. This is because all of the possible strategies are within the feasible solution set of the model and, thus, once solved, the solution defines the optimal strategy. This showcases both the benefits of posing such decision problems as MIP formulations and the range of applications that can be tackled by decision programming. Moreover, it also highlights the importance of our methodological contributions to decision programming as enablers to modelling and solving more complex, realistic problems. 

This paper is structured as follows. Section \ref{sec:decision_programming} presents the technical details associated with the decision programming framework. In Section \ref{sec:formulations}, we present the novel formulation proposed in this paper, followed by the proposed adaptation of the single policy update heuristic presented in Section \ref{sec:heuristics}.
In Section \ref{sec:computational_experiments}, we provide computational results showcasing the benefits of the methodological innovations proposed in this paper and in Section \ref{sec:case_study} we describe the case study considering optimal preventive care strategies for CHD. Lastly, in Section \ref{sec:conclusions} we provide conclusions and discuss some potential directions for further research.

\section{Decision Programming} \label{sec:decision_programming}

Decision programming relies on influence diagrams, which are graphical representations of decision problems. In influence diagrams, nodes represent chance events, decisions and consequences. Specifically, let $G(N,A)$ be an acyclic graph formed by nodes in $N = C \cup D \cup V$, where $C$ is a subset of chance nodes, $D$ a subset of decision nodes, and $V$ a subset of value nodes. Value nodes represent consequences incurred from decisions made at nodes $D$ and chance events observed at nodes $C$. Each decision and chance node $j \in C \cup D$ can assume a state $s_j$ from a discrete and finite set of states $S_j$. For a decision node $j \in D$, $S_j$ represents the decision alternatives. For a chance node $j \in C$, $S_j$ is the set of possible outcomes. 

In the diagram, arcs represent interdependency among decisions and chance events. Set $A = \{(i,j) \mid i,j \in N\}$ contains the arcs $(i,j)$, which represent the influence between nodes $i$ and $j$. This influence is propagated in the diagram in the form of \emph{information}. That is, an arc $(i,j)$ that points to a decision node $j \in D$ indicates that the decision at $j \in D$ is made \emph{knowing} the realisation (i.e., uncertainty outcome or decision made) of state $s_i \in S_i$, with $i \in C \cup D$. On the other hand, an arc that points to a chance node $j \in C$ indicates that the realisation $s_j \in S_j$ is dependent (or conditional) on realisation $s_i \in S_i$ of node $i \in C \cup D$.

The \emph{information set} $I(j) = \{i \in N \mid (i,j) \in A\}$ comprises all immediate predecessors (or 
parents) 
of a given node $j \in N$. Despite being a less common terminology, we opt for the term ``information set'' to highlight the role of information in the modelling of the decision process. The decisions $s_j \in S_j$ made in each decision node $j \in D$ depend on their \emph{information state} $s_{I(j)} \in S_{I(j)}$, where $S_{I(j)} = \prod_{i \in I(j)} S_i$ is the set of all possible information states for node $j$. Analogously, the possible realisations $s_j \in S_j$ for each chance node $j \in C$ and their associated probabilities also depend on their information state $s_{I(j)} \in S_{I(j)}$.

Let us define $X_j \in S_j$ as the realised state at a chance node $j \in C$. For a decision node $j \in D$, let $Z_j: S_{I(j)} \to S_j$ be a mapping between each information state $s_{I(j)} \in S_{I(j)}$ and decision $s_j \in S_j$. That is, $Z_j(s_{I(j)})$ defines a local decision strategy, which represents the choice of some $s_j \in S_j $ in $j \in D$, given the information $s_{I(j)}$. Such a mapping can be represented by an indicator function $\mathbb{I}: S_{I(j)} \times S_j \to \{0,1\}$ defined so that
\begin{align*}
    \mathbb{I}(s_{I(j)}, s_j) = \begin{cases} 1, &\text{ if } Z_j \text{ maps } s_{I(j)} \text{ to } s_j \text{, i.e., } Z_j(s_{I(j)}) = s_j; \\ 0, &\text{ otherwise.} \end{cases}
\end{align*}
A (global) \emph{decision strategy} is the collection of local decision strategies in all decision nodes: $Z = (Z_j)_{j \in D}$, selected from the set of all possible strategies $\mathbb{Z}$.

A \emph{path} is a sequence of states $s = (s_i)_{i=1,\dots,n}$, with $n = |C| + |D|$ and 
\begin{equation}
    S = \{(s_i)_{i=1,\dots,n} \mid s_i \in S_i, i =1, \dots, n\}\label{eq:paths}    
\end{equation} 
is the set of all possible paths. We assume that the nodes $C \cup D$ are numbered from 1 to $n$ such that for each arc $(i,j) \in A$, $i<j$. Moreover, we say that a strategy $Z$ is compatible with a path $s \in S$ if $Z_j(s_{I(j)}) = s_j$ for all $j \in D$. We denote as $S(Z) \subseteq S$ the subset of all paths that are compatible with a strategy $Z$. 

Using the notion of information states, the conditional probability of observing a given state $s_j$ for $j \in C$ is $\mathbb{P}(X_j = s_j \mid X_{I(j)} = s_{I(j)})$. The probability associated with a path $s \in S$ being observed given a strategy $Z$ can then be expressed as
\begin{align}
    \mathbb{P}(s \mid Z) = \left(\prod_{j \in C}\mathbb{P}(X_j = s_j \mid X_{I(j)} = s_{I(j)})\right)\left(\prod_{j \in D} \mathbb{I}(s_{I(j)},s_j)\right) \label{eq:path_probability}.
\end{align}
Notice that the term $\prod_{j \in D} \mathbb{I}(s_{I(j)},s_j)$ in equation \eqref{eq:path_probability} takes value one if the strategy $Z$ is compatible with the path $s \in S$, being zero otherwise. Furthermore, notice that one can pre-calculate the probability 
\begin{equation}
    p(s) = \left(\prod_{j \in C}\mathbb{P}(X_j = s_j \mid X_{I(j)} = s_{I(j)})\right) \label{eq:p-def}
\end{equation}
of a path $s \in S$ being observed, in case a compatible strategy is chosen.

At the value node $v \in V$, a real-valued utility function $U_v : S_{I(v)} \to \mathbb{R}$ maps the information state $s_{I(v)}$ to a utility value $U_v(s_{I(v)})$. We usually assume the utility value of a path $s$ to be the sum of individual value nodes' utilities: $U(s) = \sum_{v \in V} U_v(s_{I(v)})$. The default objective is to choose a strategy $Z \in \mathbb{Z}$ maximising the expected utility, which can be expressed as
\begin{equation} 
\underset{Z \in \mathbb{Z}}{\text{max }} \sum_{s \in S}  \mathbb{P}(s \mid Z) U(s). \label{eq:orig-obj}
\end{equation}
Notice that other objective functions can also be modelled. For example, \citet{salo2022} discuss the use of the conditional value-at-risk.

To formulate this into a mathematical optimisation problem, we start by representing the local strategies $Z_j$ using binary variables $z(s_j \mid s_{I(j)})$ that take value one if $\mathbb{I}(s_{I(j)},s_j) = 1$, and 0 otherwise. We then observe that using \eqref{eq:path_probability} and \eqref{eq:p-def}, the objective function \eqref{eq:orig-obj} becomes 
\begin{equation*} 
\underset{z}{\text{max }} \sum_{s \in S}  p(s) U(s) \prod_{j \in D} z(s_j \mid s_{I(j)}). \label{eq:nonlin-obj}
\end{equation*}
This function is nonlinear and is used only for illustrating the nature of the formulations. The usefulness of this construction becomes more obvious in Section \ref{sec:formulations}. \citet{salo2022} instead replace the conditional path probability $\mathbb{P}(s \mid Z)$ in \eqref{eq:orig-obj} with a continuous decision variable $\pi(s)$, enforcing the correct behaviour of this variable using affine constraints.

With these building blocks, the problem can be formulated as a mixed-integer linear programming (MILP) model, which allows for employing off-the-shelf mathematical programming solvers. The MILP problem presented in \citet{salo2022} can be stated as \eqref{eq:dp_obj}-\eqref{eq:dp_z_bin}.
\begin{align}
    \underset{Z \in \mathbb{Z}}{\text{max }}  &\sum_{s \in S} \pi(s) U(s)\label{eq:dp_obj}\\
    \text{subject to }  &\sum_{s_j \in S_j} z(s_j \mid s_{I(j)}) = 1, &&\forall j \in D, s_{I(j)} \in S_{I(j)}, \label{eq:dp_z_sum}\\
    &0 \le \pi(s) \le p(s), &&\forall s \in S, \label{eq:dp_pi_lim}\\
    &\pi(s) \le z(s_j \mid s_{I(j)}), &&\forall j \in D, s \in S, \label{eq:dp_pi_upper}\\
    &\pi(s) \ge p(s) + \sum_{j \in D} z(s_j \mid s_{I(j)}) \, - | D |, &&\forall s \in S, \label{eq:dp_pi_lower} \\
    &z(s_j \mid s_{I(j)}) \in \{0,1\}, &&\forall j \in D, s_j \in S_j, s_{I(j)} \in S_{I(j)}. \label{eq:dp_z_bin}
\end{align}

Variables $\pi(s)$ are nonnegative continuous variables representing the conditional path probability in equation \eqref{eq:path_probability}. They take the value of the path probability $p(s)$ in case the selected strategy $Z$ is compatible with the path $s \in S$ and zero otherwise. Notice that this compatibility is equivalent to observing $z(s_j \mid s_{I(j)}) = 1$ for all $s_j \in S$ such that $j \in D$. 

The objective function \eqref{eq:dp_obj} defines the expected utility value, which is calculated considering only the paths that are compatible with the strategy. Constraint \eqref{eq:dp_z_sum} enforces the one-to-one nature of the mapping $\mathbb{I}(s_{I(j)}, s_j)$, represented by the $z$-variables. The correct behaviour of variables $\pi(s)$ is guaranteed by constraints \eqref{eq:dp_pi_lim}-\eqref{eq:dp_pi_lower}, which enforce that $\pi(s) = p(s)$ if $z(s_j \mid s_{I(j)}) = 1$ for all $s_j \in S$ such that $j \in D$. The term $| D |$ in \eqref{eq:dp_pi_lower} represents the cardinality of the set $D$, that is, the number of decision nodes in the diagram. Notice that the domain of $\pi(s)$ is defined in \eqref{eq:dp_pi_lim}.

\section{Improved formulations}
\label{sec:formulations}

One key challenge associated with formulation \eqref{eq:dp_obj}--\eqref{eq:dp_z_bin}, and, in fact, any MILP formulation, is that computational performance is strongly tied to the tightness of the formulation. In this context, the tightness of a MILP formulation is related to how close the linear relaxation solution is to the initial primal bound, e.g., the first integer feasible solution value obtained by the solver during the solution process or one obtained using primal heuristics. 

Next, we present reformulations developed to enhance the numerical performance of the decision programming formulation \eqref{eq:dp_obj}--\eqref{eq:dp_z_bin}. For that, let us first define the subset of paths
\begin{equation*}
    S_{s_j \mid s_{I(j)}} = \braces{s \in S \mid (s_{I(j)}, s_j) \subseteq s}.    
\end{equation*}
Notice that we use the notation $(s_{I(j)}, s_j)$ to represent a portion of a path $s$, formed by the combination of the information state $s_{I(j)}$ (which may itself be a collection of states, if $|I(j)| > 1$) and the state $s_j$. We also utilise the set operator $\subseteq$ to indicate that the states $(s_{I(j)}, s_j)$ are part of the path $s \in S$. Notice that the states $(s_{I(j)}, s_j)$ do not need to be consecutive in the path $s$, although the ordering between $s_{I(j)}$ and $s_j$ is naturally preserved in $s$. 

Considering $j \in D$, the subset $S_{s_j \mid s_{I(j)}}$ allows us to define the notion of \emph{locally compatible paths}, that is, the collection of paths $s$ compatible with local strategies $Z_j$ for which $\mathbb{I}(s_{I(j)}, s_j) = 1$. The definition of the subset $S_{s_j \mid s_{I(j)}}$ allows us to derive the following valid inequality for \eqref{eq:dp_obj}--\eqref{eq:dp_z_bin}: 
\begin{equation} \label{eq:locally_compatible_paths}
    \sum_{s \in S_{s_j \mid s_{I(j)}}} \pi(s) \leq z(s_j \mid s_{I(j)}), \quad \forall j \in D, s_j \in S_j, s_{I(j)} \in S_{I(j)}.
\end{equation}
Constraint \eqref{eq:locally_compatible_paths} states that only paths that are compatible with the selected strategy might be allowed to have a probability different than zero. Moreover, since it is enforced on all decision nodes, it means that this constraint guarantees that only the paths that are compatible with the strategy $Z$ are active. Recall that we denote this set of compatible paths as $S(Z) \subseteq S$.

As pointed out in \citet{salo2022}, for expected utility maximization, constraint \eqref{eq:dp_pi_lower}, which prevents variables $\pi(s)$ from wrongly taking value zero, is only required when some of the utility values $U(s)$, $s \in S$, are negative. Notice that this is otherwise prevented by the maximization of the objective function \eqref{eq:dp_obj}, naturally steering these variables to their upper bound values. Another way to guarantee that the variables $\pi(s)$ take their correct value, i.e., $\pi(s) = p(s)$, if $s \in S(Z)$, is to impose the constraint
\begin{equation} \label{eq:total_prob}
    \sum_{s \in S} \pi(s) = 1.
\end{equation}

As it will be discussed in Section \ref{sec:computational_experiments}, replacing \eqref{eq:dp_pi_upper} and \eqref{eq:dp_pi_lower} with \eqref{eq:locally_compatible_paths} and \eqref{eq:total_prob} provides considerable gains in terms of linear relaxation strengthening. Furthermore, we observe that the computational performance can be even further improved by employing a simple variable substitution. Recall that in the original formulation \eqref{eq:dp_obj}-\eqref{eq:dp_z_bin}, variables $\pi(s)$ represent the conditional path probability $\mathbb{P}(s \mid Z) = p(s) \prod_{j \in D} z(s_j \mid s_{I(j)})$. If we let $x(s) \in [0,1]$, $s \in S$ represent the product $\prod_{j \in D} z(s_j \mid s_{I(j)})$, then we can reformulate the problem by substituting $\pi(s) = p(s)x(s)$ for all $s \in S$.

Although $x(s)$, $s \in S$, is continuous, it behaves as a binary variable which takes value one whenever the path is compatible with the strategy and zero, otherwise. This is analogous to the behaviour of variable $\pi(s) \in [0,p(s)]$ in \eqref{eq:dp_obj}--\eqref{eq:dp_z_bin}. We highlight that, from a theoretical standpoint, there is no obvious reason for performing such a substitution. On the other hand, we will show in the computational experiments presented in Section \ref{sec:case_study} that it yields significant practical benefits in terms of computational performance.

Using these $x$-variables, we can reformulate \eqref{eq:locally_compatible_paths} as
\begin{equation} \label{eq:dp2_locally_compatible_paths_1}
    \sum_{s \in S_{s_j \mid s_{I(j)}}} x(s) \leq |S_{s_j \mid s_{I(j)}}|z(s_j \mid s_{I(j)}), \quad \forall j \in D, s_j \in S_j, s_{I(j)} \in S_{I(j)}, 
\end{equation}
a consequence of $x(s) \in [0,1]$ and the fact that $z(s_j \mid s_{I(j)})$ must be equal to 1 for $x(s)$ to be positive for $s \in S_{s_j \mid s_{I(j)}}$.

Constraint \eqref{eq:dp2_locally_compatible_paths_1} can be strengthened further. We note that a path must be in the set of compatible paths $S(Z)$ in order for $x(s)$ to be positive with strategy $Z$. Using this information, we can infer a tighter upper bound for the number of paths that can be active ($x(s) >0$) from the set of locally compatible paths. We observe that in a set of compatible paths $S(Z)$, each information state $s_{I(j)}$ maps to exactly one decision alternative $s_j$ for each decision node $j \in D$, in accordance with constraint \eqref{eq:dp_z_sum}. However, the set of locally compatible paths for a given pair of information state and decision node state $(s_{I(j)}, s_j)$ of decision node $j \in D$, includes paths for all combinations $(s_{I(k)}, s_k)$ of information states and decisions for the other decision nodes $k \in D \setminus \braces{j}$. Hence, only a fraction of the locally compatible paths can be active. The fraction is linked to the number of states $| S_k |$ of the other decision nodes $k \in D \setminus \braces{j}$. The number of locally compatible paths that will also be active, i.e., $|S_{s_j \mid s_{I(j)}} \cap S(Z)|$ can be defined as
\begin{equation} \label{eq:dp2_locally_active_compatible_paths}
    |S_{s_j \mid s_{I(j)}} \cap S(Z)| = \frac{|S_{s_j | s_I(j)}|} {\Pi_{k \in D \setminus (\braces{j} \cup I(j))} |S_k|}.
\end{equation}
Notice that the calculation of the number of active paths must take into account the fact that some decision nodes may be part of the information state $I(j)$ of node $j \in D$, and, as such, will have their states observed (or fixed) in the set $S_{s_j | s_I(j)}$. Therefore, these decision nodes must be excluded from the product in the denominator in equation \eqref{eq:dp2_locally_active_compatible_paths}. Using \eqref{eq:dp2_locally_active_compatible_paths}, we can reformulate \eqref{eq:dp2_locally_compatible_paths_1} into the strengthened form
\begin{equation} \label{eq:dp2_locally_compatible_paths_2}
    \sum_{s \in S_{s_j \mid s_{I(j)}}} x(s) \leq  \frac{|S_{s_j | s_I(j)}|} {\Pi_{k \in D \setminus (\braces{j} \cup I(j))} |S_k|}z(s_j \mid s_{I(j)}), \quad \forall j \in D, s_j \in S_j, s_{I(j)} \in S_{I(j)}. 
\end{equation}

One last aspect that can be taken into account is that, depending on the problem structure, some sequence of states $s = (s_i)_{i=1,\dots,n}$ forming a path may never be observed and can be preemptively filtered out from the set of paths $S$. This is the case, for example, in problems where earlier decisions or uncertain events dictate whether alternatives or uncertainties are observed. For instance, an initial decision regarding whether or not to build an industrial plant naturally restricts subsequent decisions regarding capacity expansion. Analogously, it may be that an uncertain production rate is only observed if one decides to build the production facility in the first place. To prevent the assembling of these unnecessary paths, we consider a set of \emph{forbidden} paths, which, once removed, lead to a set $S^* \subseteq S$ of \emph{effective} paths. Notice that these forbidden paths have probability zero by the structure of the problem, and therefore their removal does not affect the expected utility nor the constraints of the model. Furthermore, their removal allows for significant savings in terms of the scale of the model.

One issue emerges in settings where $S^* \subset S$ regarding the term \eqref{eq:dp2_locally_active_compatible_paths}. Notice that the bound is based on the premise that we can infer the total number of paths by considering the Cartesian product of the state sets $S_j$, $j \in N$. However, as forbidden paths are removed, some of the $x$-variables corresponding to paths $s \in S_{s_j | s_I(j)}$ might be removed, making inequality \eqref{eq:dp2_locally_compatible_paths_2} loose. A simple safeguard for this is to consider  
\begin{equation} \label{eq:gamma_value}
    \Gamma(s_j|s_{I(j)}) = \min \braces{|S^*_{s_j \mid s_{I(j)}}|, \frac{|S_{s_j | s_I(j)}|} {\Pi_{k \in D \setminus (\braces{j} \cup I(j))} |S_k|}}
\end{equation}
and reformulate \eqref{eq:dp2_locally_compatible_paths_2} as
\begin{equation} \label{eq:dp2_locally_compatible_paths_3}
    \sum_{s \in S_{s_j \mid s_{I(j)}}} x(s) \leq \Gamma(s_j|s_{I(j)})z(s_j \mid s_{I(j)}), \quad \forall j \in D, s_j \in S_j, s_{I(j)} \in S_{I(j)}. 
\end{equation}

Combining the above, we can reformulate \eqref{eq:dp_obj}--\eqref{eq:dp_z_bin} as follows.
\begin{align}
    \underset{Z \in \mathbb{Z}}{\text{maximize }}  &\sum_{s \in S^*} U(s)p(s)x(s)\label{eq:dp2_obj}\\
    \text{subject to }  
    &\sum_{s_j \in S_j} z(s_j \mid s_{I(j)}) = 1, &&\forall j \in D, s_{I(j)} \in S_{I(j)} \label{eq:dp2_z_sum}\\
    & \sum_{s \in S_{s_j \mid s_{I(j)}}} x(s) \leq \Gamma(s_j|s_{I(j)}) z(s_j \mid s_{I(j)}), &&\forall j \in D, s_j \in S_j, s_{I(j)} \in S_{I(j)} \label{eq:dp2_x_upper} \\
    & \sum_{s \in S^*} p(s)x(s) = 1, && \label{eq:dp2_prob_sum} \\
    & 0 \leq x(s) \leq 1, &&\forall s \in S^* \label{eq:dp2_x_lim} \\ 
    &z(s_j \mid s_{I(j)}) \in \{0,1\}, &&\forall j \in D, s_j \in S_j, s_{I(j)} \in S_{I(j)} \label{eq:dp2_z_bin}.
\end{align}
where $\Gamma(s_j|s_{I(j)})$ is defined as in \eqref{eq:gamma_value}. Note that this formulation preserves the (mixed-integer) linear nature of \eqref{eq:dp_obj}--\eqref{eq:dp_z_bin}.

As discussed earlier, one of the main advantages of the decision programming formulation is the ability to incorporate objectives and constraints involving arbitrary utility functions and probability distributions within the model. \citet{salo2022} demonstrate this by proposing a model that considers conditional value-at-risk as one of the utility functions. The same can be achieved with our proposed formulation, by simply substituting $p(s)x(s)$ in place of variables $\pi(s)$. 

The path-based structure of \eqref{eq:dp_obj}--\eqref{eq:dp_z_bin} makes formulating chance and budget constraints straightforward. For modeling chance constraints, we can define $\tilde{S}$ as the set of ``undesirable'' paths, the total probability of which must not exceed $\rho$. The corresponding chance constraint is then 
\begin{equation} \label{eq:chance_constraint}
    \sum_{s \in \tilde{S}} x(s)p(s) \le \rho. 
\end{equation}

Likewise, the set $\tilde{S}$ could further be defined as, e.g., the set of paths with a small utility $U(s) \le u_{threshold}$. With $\rho=0$, \eqref{eq:chance_constraint} can be seen as a budget constraint, stating that for all compatible paths $s \in S(Z)$ with $p(s)>0$, the utility $U(s)$ must be at least $u_{threshold}$.



\section{Primal heuristic: single policy update (SPU)} 
\label{sec:heuristics}

While our approach of formulating influence diagrams into mixed-integer linear models does allow us to use powerful off-the-shelf solvers, it is still hindered by the well-known fact that solving such problems is NP-hard \citep{schrijver2003combinatorial}. To make MIP solvers more efficient, \emph{primal heuristics} are used to obtain and improve integer solutions. Obtaining good starting integer solutions can have a significant impact on the performance of branch-and-bound solvers, as it helps in pruning poor-quality solutions early.

Decision programming is based on limited-memory influence diagrams (LIMIDs) and solution approaches presented in previous literature can be used to obtain solutions to these problems. A notable contribution of \citet{lauritzen2001representing} is the single policy update (SPU) heuristic for obtaining ``locally optimal'' strategies in the sense that the corresponding solutions cannot be improved by changing only one of the local strategies $Z_j(s_{I(j)})$.  

Our proposed heuristic is loosely based on the ideas in \citet{lauritzen2001representing}, as described in Algorithm \ref{alg:spu}. The first step of the heuristic is to obtain a random strategy $Z$ (note that this too is a heuristic, albeit a very simple one). Additionally, we initialise the $lastImprovement$ variable that will be used to stop the algorithm after finding a local optimum. The strategy $Z$ is then iteratively improved by examining each information state $s_{I(j)} \in S_{I(j)}$ for each decision node $j \in D$ in order, choosing the local strategy $Z'_j(s_{I(j)})$ maximising the expected utility. We obtain incrementally improving strategies by replacing the local strategy $Z_j(s_{I(j)})$ with $Z'_j(s_{I(j)})$ whenever the change results in an increase in expected utility. Finally, the pair $(j, s_{I(j)})$ is stored in the $lastImprovement$ variable if an improvement has been observed.

\begin{algorithm}
\caption{The single policy update heuristic}\label{alg:spu}
$Z \gets randomstrategy()$\;
$lastImprovement \gets (undef, undef)$\;
\While{true }{
    \For{$j \in D$, $s_{I(j)} \in S_{I(j)}$}{
        \eIf{$(j,s_{I(j)}) = lastImprovement$}{
            \KwRet{$Z$}\;
        }{
            $Z'_j(s_{I(j)}) \gets bestLocalStrategy(Z,j,s_{I(j)})$\;
            $Z' \gets modifyStrategy(Z,Z'_j(s_{I(j)}))$\;
            \If{$EU(Z') > EU(Z)$}{
                $Z \gets Z'$\;
                $lastImprovement \gets (j, s_{I(j)})$\;
            }
        }
    }
}
\end{algorithm}

This process of locally improving the strategy is performed repeatedly for all pairs $(j, s_{I(j)})$ until no improvement is made during a whole iteration through the set of such pairs, that is, $(j, s_{I(j)}) = lastImprovement$. The number of possible strategies $Z$ is finite, and the algorithm thus converges in a finite number of iterations. It is also easy to see that at termination, there is no possible local improvement and the strategy $Z$ is thus, in that sense, locally optimal. \citet{lauritzen2001representing} show that for \emph{soluble} LIMIDs, this heuristic results in a globally optimal solution. However, influence diagrams are not generally soluble. The performance of the heuristic is explored in Section \ref{sec:computational_experiments}.

\section{Computational experiments} \label{sec:computational_experiments}

We present a collection of computational experiments carried out to assess the performance of the proposed reformulation \eqref{eq:dp2_obj}-\eqref{eq:dp2_z_bin} against the original formulation \eqref{eq:dp_obj}-\eqref{eq:dp_z_bin} presented in \citet{salo2022}. 
In addition, we present computational results highlighting the performance of the proposed SPU heuristic. The problems used for testing are (i) the pig farm problem originally from \citet{lauritzen2001representing} modified to allow for artificially augmenting the number of time periods and generating input parameters randomly; (ii) the N-monitoring problem, as proposed in \citet{salo2022}, in which we also can artificially augment the number of decision agents and randomly generate instances. These two problems have a major difference in their structure: while the N-monitoring problem has $N$ decisions in parallel with no communication between the decision-makers, the pig farm problem is a partially observed Markov decision process (POMDP) where decisions are made in series with limited memory of the past. This structural difference leads to the N-monitoring problem having a larger treewidth, which has generally been an issue for solving influence diagrams \citep{maua2012solving}.
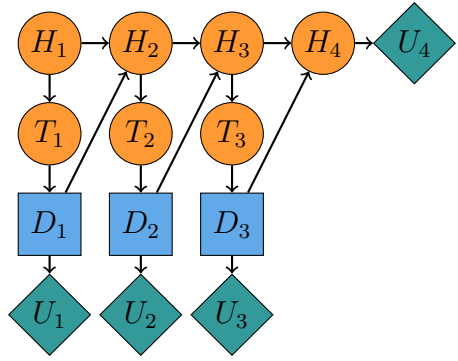
\begin{figure}
    \centering
    \begin{tikzpicture}
    [decision/.style={fill=blue!80, draw, minimum size=2em, inner sep=2pt}, 
    chance/.style={circle, fill=orange!80, draw, minimum size=2em, inner sep=2pt},
    value/.style={diamond, fill=teal!80, draw, minimum size=2em, inner sep=2pt},
    scale=1.2]
    \node[chance]   (H1) at (0, 5)  {$H_1$};
    \node[chance]   (T1) at (0, 4)  {$T_1$};
    \node[decision] (D1) at (0, 3)  {$D_1$};
    \node[value]    (U1) at (0, 2)  {$U_1$};
    
    \node[chance]   (H2) at (1, 5)  {$H_2$};
    \node[chance]   (T2) at (1, 4)  {$T_2$};
    \node[decision] (D2) at (1, 3)  {$D_2$};
    \node[value]    (U2) at (1, 2)  {$U_2$}; 
    
    \node[chance]   (H3) at (2, 5)  {$H_3$};
    \node[chance]   (T3) at (2, 4)  {$T_3$};
    \node[decision] (D3) at (2, 3)  {$D_3$};
    \node[value]    (U3) at (2, 2)  {$U_3$}; 
    
    \node[chance]   (H4) at (3, 5)  {$H_4$};
    \node[value]    (U4) at (4, 5)  {$U_4$}; 
    
    \draw[->, thick] (H1) -- (T1);
    \draw[->, thick] (H1) -- (H2);
    \draw[->, thick] (T1) -- (D1);
    \draw[->, thick] (D1) -- (U1);
    \draw[->, thick] (D1) -- (H2);
    
    \draw[->, thick] (H2) -- (T2);
    \draw[->, thick] (H2) -- (H3);
    \draw[->, thick] (T2) -- (D2);
    \draw[->, thick] (D2) -- (U2);
    \draw[->, thick] (D2) -- (H3);
    
    \draw[->, thick] (H3) -- (T3);
    \draw[->, thick] (H3) -- (H4);
    \draw[->, thick] (T3) -- (D3);
    \draw[->, thick] (D3) -- (U3);
    \draw[->, thick] (D3) -- (H4);

    \draw[->, thick] (H4) -- (U4);
\end{tikzpicture}
    \caption{The influence diagram of the pig farm problem \cite{salo2022}.}
    \label{fig:pigfarm}
\end{figure}

The pig farm problem is presented in Figure \ref{fig:pigfarm}, where nodes $H_i$ represent the health of a pig, nodes $T_i$ represent the result of a diagnostic test, and nodes $D_i$ correspond to treatment decisions. Treating a pig incurs a cost represented by the value nodes $U_i$ for $i \in \{1,2,3\}$ and in the end, the pig can be sold for a price depending on the final health state. This problem was chosen for the computational comparison because similar low treewidth structures frequently arise in contexts such as quality control \citep{cobb2024intermittent} or testing and treating patients for a disease \citep{hynninen2019value}, discussed in Section \ref{sec:case_study} of this paper. 

All experiments are run using 16GB of memory and 8 threads on an Intel Xeon Gold 6248 CPU and the code can be found in \\ \href{https://github.com/gamma-opt/DecisionProgramming.jl/tree/new-formulation-experiments}{https://github.com/gamma-opt/DecisionProgramming.jl/tree/new-formulation-experiments}.

\subsection{Problem size}

First, we compare the model sizes of the two formulations presented in Sections \ref{sec:decision_programming} and \ref{sec:formulations}. In both formulations, the number of variables is the same. There are $\sum_{j \in D} |S_j||S_{I(j)}|$ $z$-variables and $|S|$ path variables, either $\pi$ or $x$, depending on the formulation. As for the number of constraints, the formulation \eqref{eq:dp_obj}-\eqref{eq:dp_z_bin} has $\sum_{j \in D}|S_{I(j)}|$ constraints \eqref{eq:dp_z_sum}, $2|S|$ bounds for $\pi$-variables, $|D||S|$ constraints \eqref{eq:dp_pi_upper} and $|S|$ constraints \eqref{eq:dp_pi_lower}. Arranging the terms, the total number of constraints becomes 
\begin{equation}
    \label{eq:num_paths}
    (3+|D|)|S| + \sum_{j \in D}|S_{I(j)}|.
\end{equation}

The formulation \eqref{eq:dp2_obj}-\eqref{eq:dp2_z_bin} has $\sum_{j \in D}|S_{I(j)}|$ constraints \eqref{eq:dp2_z_sum}, $\sum_{j \in D}|S_j||S_{I(j)}|$ constraints \eqref{eq:dp2_x_upper}, one constraint \eqref{eq:dp2_prob_sum} and $2|S|$ bounds for $x$-variables. Arranging the terms, the total number of constraints becomes 
\begin{equation}
    \label{eq:num_paths_2}
    2|S| + \sum_{j \in D}(1+|S_j|)|S_{I(j)}|.
\end{equation}

We note that $|S|=\prod_{j \in C \cup D}|S_j|$ and that especially with a large number of nodes, the first term becomes impractically large in both \eqref{eq:num_paths} and \eqref{eq:num_paths_2}. The increase in the number of path-related constraints is exponential, while the increase in the rest of the constraints is often linear, as shown in the following two example problems.

\subsubsection{Pig farm}

For the pig farm example presented in \citet{lauritzen2001representing}, we observe that the problem consists of $3n+1$ decision and chance nodes, where $n$ is the number of decision stages, and that $|S_j|=2$ for all nodes $j \in C \cup D$. Note that this is slightly different to \citet{lauritzen2001representing}, where the length of the problem is tied to the number of health nodes. The length of a problem with $n$ decision nodes would then be $n+1$.

With these observations, $|S|=2^{3n+1}$ and $|D|=n$. Thus, the number of constraints in Eq. \eqref{eq:num_paths} becomes $(3+n)2^{3n+1} + 2n$ and the corresponding number in Eq. \eqref{eq:num_paths_2} becomes $2^{3n+2} + 6n$.

\subsubsection{N-monitoring}

 We can perform a similar analysis for the N-monitoring example presented in \citet{salo2022}. The problem consists of $2n+2$ decision and chance nodes, where $n$ is the number of report-action pairs. As in the pig farm problem, $|S_j|=2$ for all nodes $j \in C \cup D$. 

With these observations, $|S|=2^{2n+2}$ and $|D|=n$. Thus, the number of constraints in Eq. \eqref{eq:num_paths} becomes $(3+n)2^{2n+2} + 2n$ and the corresponding number in Eq. \eqref{eq:num_paths_2} becomes $2^{2n+3} + 6n$.  

\subsection{Solution times}

\begin{figure}[ht]
\centering
     \begin{subfigure}[b]{0.48\textwidth}
         \centering
         \includegraphics[width=\textwidth]{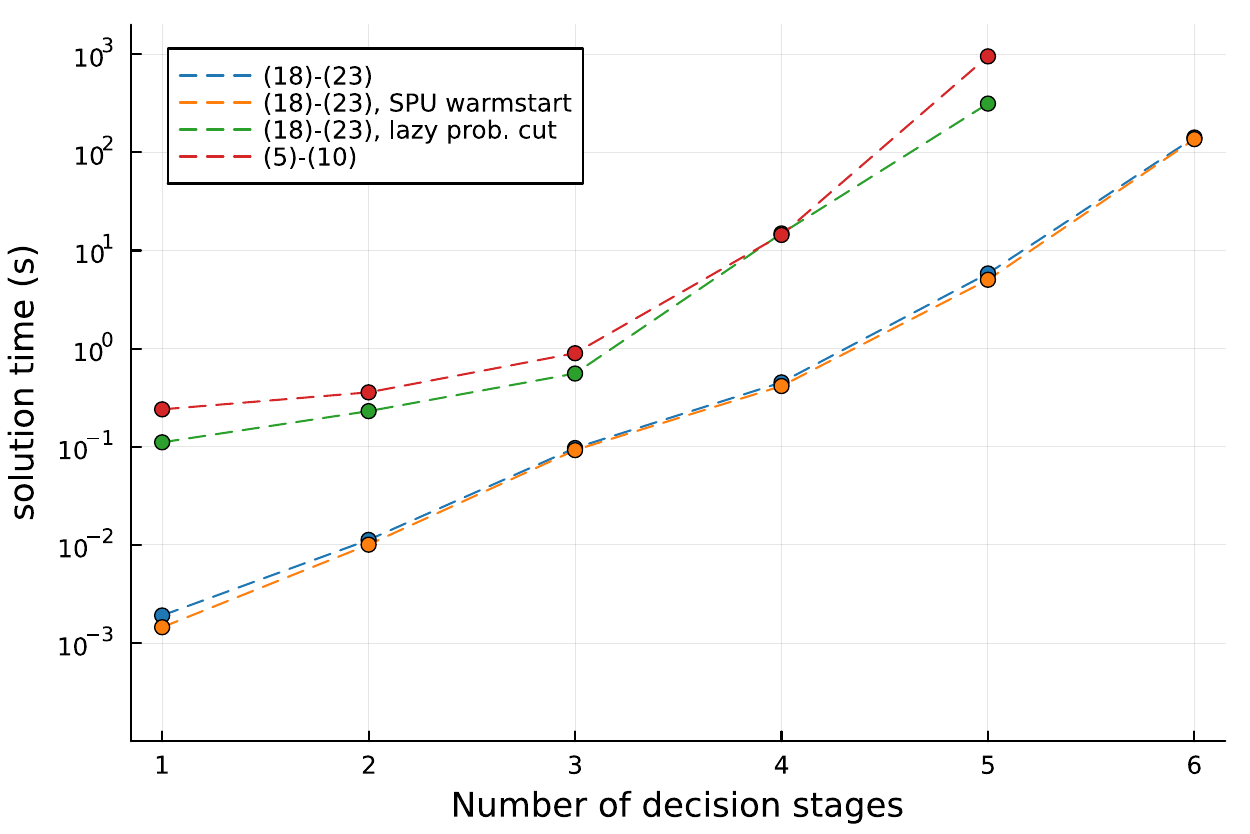}
         \caption{The pig farm problem}
         \label{fig:sol_times_pigfarm}
     \end{subfigure}
     \hfill
     \begin{subfigure}[b]{0.48\textwidth}
         \centering
         \includegraphics[width=\textwidth]{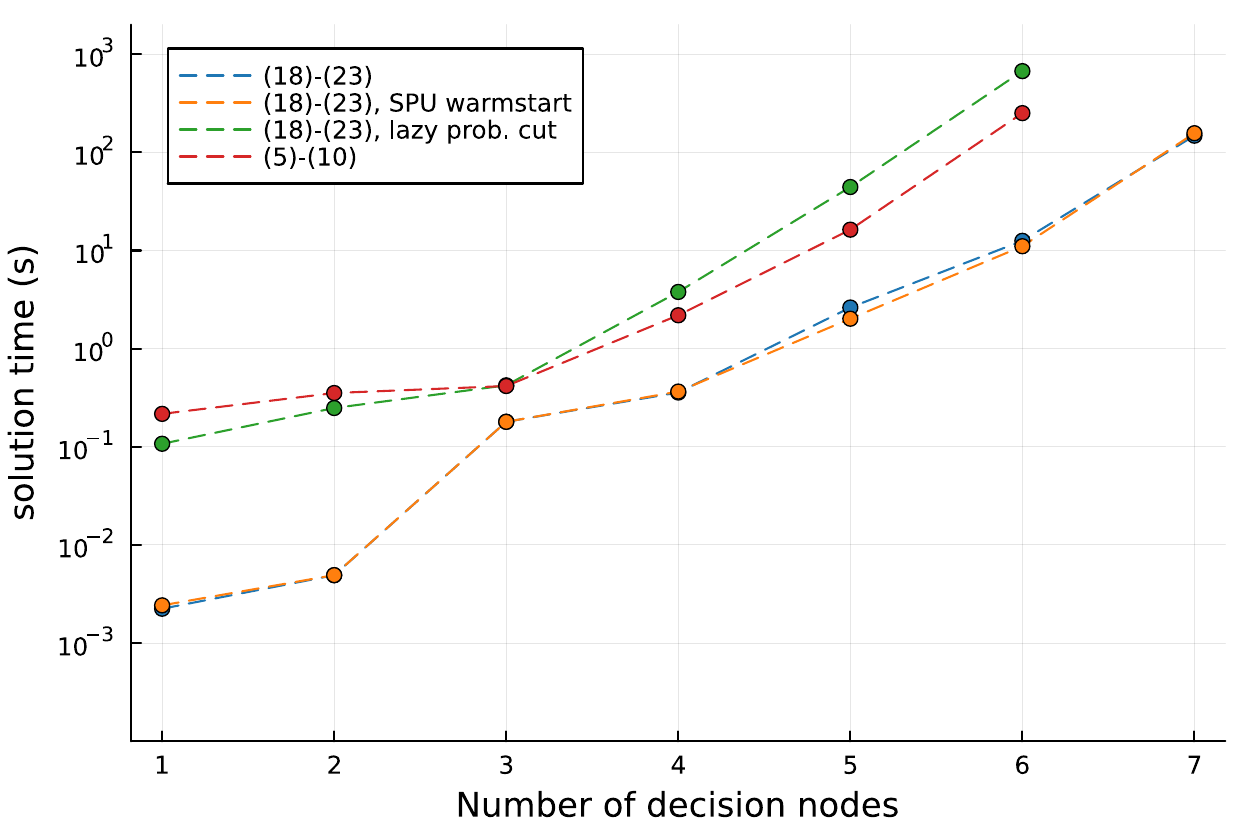}
         \caption{The N-monitoring problem}
         \label{fig:sol_times_nmonitoring}
     \end{subfigure}
     \caption{The solution times of the two example problems with different number of decision nodes using different formulations. Notice the logarithmic y-axis.}
     \label{fig:sol_times}
\end{figure}

Figure \ref{fig:sol_times} shows the increase in average solution times over 50 instances as the number of decision stages increases in the two example problems. For the original formulation \eqref{eq:dp_obj}-\eqref{eq:dp_z_bin}, \citet{salo2022} show that solution times are greatly improved by adding a \emph{probability cut} $\sum_{s \in S} \pi(s) = 1$ as a lazy constraint to the model. A lazy constraint is a constraint that is added to the model formulation when it is deemed violated by an incumbent feasible solution found in the branch-and-cut tree search, instead of adding it from the beginning of the solution process. This approach is thus used in the computational experiments for the original formulation. For \eqref{eq:dp2_obj}-\eqref{eq:dp2_z_bin}, a similar constraint is included in the formulation by default. However, we additionally analyze an instance of the the reformulated model \eqref{eq:dp2_obj}-\eqref{eq:dp2_z_bin}, where constraint \eqref{eq:dp2_prob_sum} is implemented as a lazy constraint.

For both problems, it seems that the rate of increase in the solution times quickly renders the original formulation \eqref{eq:dp_obj}-\eqref{eq:dp_z_bin} computationally intractable, as seen in Fig \ref{fig:sol_times}. This was also noted by \citet{salo2022} in their computational results. The solution times for the improved formulation \eqref{eq:dp2_obj}-\eqref{eq:dp2_z_bin} using locally compatible path sets seem to increase at a slower rate than for the original formulation. Interestingly, the solution times are similar for both problems, suggesting that the treewidth does not have a significant impact on the computational performance, unlike for the methods discussed in \citet{lauritzen2001representing} or \citet{maua2012solving}. 
Finally, the lazy probability cut that was found to improve solution times in \citet{salo2022} is detrimental to computational performance in the new formulation \eqref{eq:dp2_obj}-\eqref{eq:dp2_z_bin}. 

In Table \ref{tbl:stats}, we present statistics on the quality of the LP relaxation. As discussed before, the hypothesis is that the formulation \eqref{eq:dp2_obj}-\eqref{eq:dp2_z_bin} is considerably tighter than \eqref{eq:dp_obj}-\eqref{eq:dp_z_bin}. The results from the pig farm problem strongly support this, as more than half of the LP relaxation solutions for the novel formulation are within 25\% of the optimal solution, while the solutions using \eqref{eq:dp_obj}-\eqref{eq:dp_z_bin} are orders of magnitude further from the optimal solution. 

\begin{table}[ht]
\centering
\begin{tabular}{l|ll}
                & \eqref{eq:dp_obj}-\eqref{eq:dp_z_bin} & \eqref{eq:dp2_obj}-\eqref{eq:dp2_z_bin} \\ \hline
10th percentile & 15.4 & 1.00 \\
median          & 26.4 & 1.21 \\
90th percentile & 31.1 & 1.81 \\
mean            & 25.0 & 1.34
\end{tabular}
\caption{Statistics of the root relaxation quality relative to the optimal solution for 50 randomly generated pig farm problems with 5 decision stages. The solutions are scaled so that a value of 1 corresponds to the optimal solution.}
\label{tbl:stats}
\end{table}

\begin{figure}[t]
     \centering
     \includegraphics[width=0.8\textwidth]{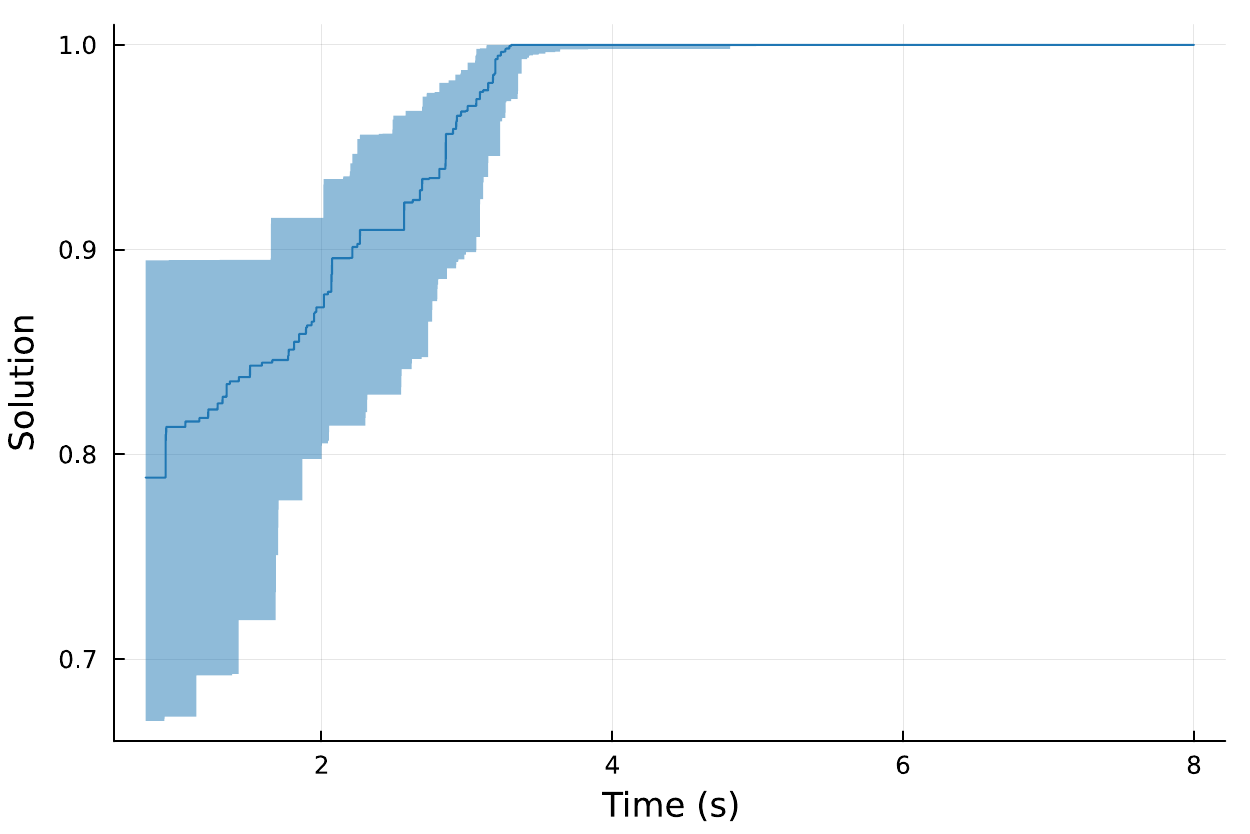}
     \caption{The median and first and third quartiles of solutions found by the SPU heuristic in 50 randomly generated pig farm problems with 6 decision stages. The solutions are scaled so that a value of 1 corresponds to the optimal solution.}
     \label{fig:spu}
 \end{figure}

Figure \ref{fig:spu} shows the process of improving solutions in the single policy update (SPU) heuristic. For the 50 instances in this test set, the last solution is found within eight seconds, and the solution is the global optimum in all but one of the instances. Note that \citet{lauritzen2001representing} showed that this version of the pig farm problem is not soluble, and thus the SPU heuristic is not guaranteed to find the optimal solution. We observe that while the single policy update heuristic is successful in finding good initial solutions quickly, the effect of providing the solver with these initial solutions is negligible (see Figure \ref{fig:sol_times}). Thus, our results suggest that improving the LP relaxation bound has a much greater impact on improving the solution time.

\section{Case study: optimal preventive healthcare for CHD} \label{sec:case_study}

\begin{figure}
    \centering
    \begin{tikzpicture}
      \draw (0,0) node[below] {\scriptsize $0\%$}; 
      \draw (5,0) node[below] {\scriptsize $100\%$};
      \filldraw[teal!80] (0, 0) rectangle +(2.25, 0.5);
      \draw (1.125,0) node[above] {\small do not treat};
      \filldraw[orange!80] (2.25, 0) rectangle +(1.25, 0.5);
      \draw (2.875,0) node[above] {\small test};
      \filldraw[blue!80] (3.5, 0) rectangle +(1.5, 0.5);  
      \draw (4.25,0) node[above] {\small treat};
      \draw[thick] (0, 0) -- (5, 0);
      \draw[thick] (0, 0.1) -- (0, -0.1) ++ (5,0) -- +(0,.2);
      \draw[thick, ->] (2.25, -0.5) node[below] {\scriptsize {$T_t$}} -- +(0, 0.5) ;
      \draw[thick, ->] (3.5, -0.5) node[below] {\scriptsize {$T_{tt}$}} -- +(0, 0.5) ;
\end{tikzpicture}
    \caption{Testing ($T_t$) and test-treatment ($T_{tt}$) thresholds.}
    \label{fig:thresholds}
\end{figure}
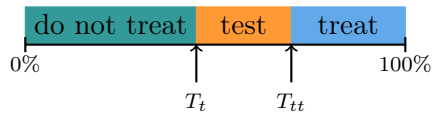

One of the first frameworks for medical decision-making considering whether to treat, test or not treat was developed by \citet{pauker1980threshold}. This framework provides an analytical basis for optimal testing and treatment strategies. They developed two thresholds, referred to as ``testing'' and ``test-treatment'' thresholds. The thresholds are probability cut-offs and they divide subjects into three groups: if the risk of disease is below the ``testing'' threshold, treatment and testing should be withheld, if it is above the ``test-treatment'' threshold, treatment should be given and if the risk falls in between these thresholds then a diagnostic test should be performed and the treatment decision made based on its results. The thresholds are visualised in Figure \ref{fig:thresholds}. 

In this case study we use decision programming to optimise the use of traditional and genetic testing to support the targeting of statin medication treatment for preventing coronary heart disease (CHD). This case study is replicated from \citet{hynninen2019value}, where the authors developed a testing and treatment strategy by optimising net monetary benefit (NMB), a cost-benefit objective consisting of the health outcomes and testing costs within a 10-year time horizon. 

The decision process stems from the patient's state of health, represented by a chance event $H$ describing whether the patient will or will not have a CHD event in the following 10 years. The probability of a CHD event is assumed to be described by a prior risk estimate $R_0$ based on factors such as the age and sex of the patient. The likelihood of a correct prognosis can be improved by carrying out tests on traditional risk factors (TRS), genetic risk factors (GRS) or both. Based on their prognosis, a decision is made on whether a patient is subjected to preventive treatment with statin medication. 

In \citet{hynninen2019value}, six predefined testing and treatment strategies were evaluated independently. In each of these strategies, the optimal allocation of tests and treatment according to risk estimates was obtained by solving the associated decision tree via dynamic programming. The six strategies considered in \citet{hynninen2019value} were: 
\begin{enumerate*}[label=(\roman*)]
    \item no tests and no treatment (‘No treatment’);
    \item using prior risk to allocate treatment (‘Treatment optimised’);
    \item performing TRS on optimised patient segment and allocating treatment based on updated risk estimates (‘TRS optimized’);
    \item performing GRS on optimised patient segment and allocating treatment based on updated risk estimates (‘GRS optimized’);
    \item performing TRS on optimised patient segment and based on its results performing GRS optimally to allocate treatment (‘TRS \& GRS optimized’);
    \item performing GRS on optimised patient segment and based on its results performing TRS optimally to allocate treatment (‘GRS \& TRS optimized’).
\end{enumerate*}

Essentially, this comprises determining optimal ``testing'' and ``test-treatment'' thresholds (cf. Figure 
\ref{fig:thresholds}) for TRS and GRS from the perspective of net monetary benefit (NMB) for each strategy (i-vi). Interestingly, the threshold values for GRS in \citet{hynninen2019value} were different than the ones found in the study presented in \citet{tikkanen2013genetic}. This is due to the different perspectives – pure patient welfare versus NMB – that the studies were conducted from. For example, the national health care guidelines for allocating treatment were not considered in the optimisation in \citet{hynninen2019value}. This showcases that the two thresholds described in \citet{pauker1980threshold} are not unique for a given disease and prognostic test because the perspective of the study affects the threshold values.

Analogously, our decision programming model determines an optimal decision strategy for allocating preventive care for CHD. The data and structure of the problem are the same as those utilised in \citet{hynninen2019value}. However, due to the flexibility of decision programming, the strategies (i-vi) do not need to be predefined. Instead, we can optimise the design of the strategy simultaneously with the threshold values, because all of these strategies are within the feasible solutions of the resulting optimisation model. 

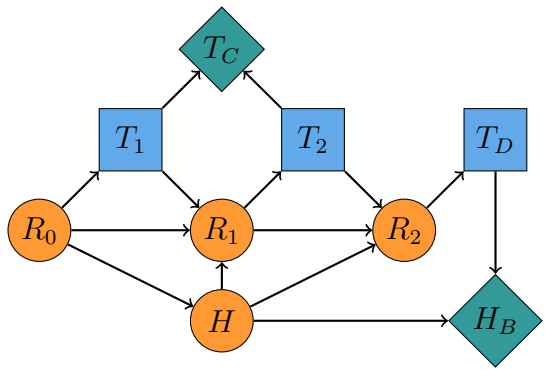
\begin{figure}[ht]
    \centering
    \begin{tikzpicture}
    [decision/.style={fill=blue!80, draw, minimum size=2em, inner sep=2pt}, 
    chance/.style={circle, fill=orange!80, draw, minimum size=2em, inner sep=2pt},
    value/.style={diamond, fill=teal!80, draw, minimum size=2em, inner sep=2pt},
    scale=1.2]
    \node[chance]   (R_0) at (0, 2)  {$R_0$};
    \node[chance]   (R_1) at (2, 2)  {$R_1$};
    \node[chance]   (R_2) at (4, 2)  {$R_2$};
    \node[chance]   (H)   at (2, 1)  {$H$};
    \node[decision] (T_D) at (5, 3)  {$T_D$};
    \node[decision] (T_1) at (1, 3)  {$T_1$};
    \node[decision] (T_2) at (3, 3)  {$T_2$};
    \node[value]    (T_C) at  (2, 4)  {$T_C$};
    \node[value]    (H_B) at  (5, 1)  {$H_B$};
    \draw[->, thick] (R_0) -- (R_1);
    \draw[->, thick] (R_0) -- (T_1);
    \draw[->, thick] (R_0) -- (H);
    \draw[->, thick] (T_1) -- (T_C);
    \draw[->, thick] (T_1) -- (R_1);
    \draw[->, thick] (R_1) -- (T_2);
    \draw[->, thick] (T_2) -- (R_2);
    \draw[->, thick] (T_2) -- (T_C);
    \draw[->, thick] (H) -- (R_1);
    \draw[->, thick] (H) -- (R_2);
    \draw[->, thick] (R_2) -- (T_D);
    \draw[->, thick] (R_1) -- (R_2);
    \draw[->, thick] (H) -- (H_B);
    \draw[->, thick] (T_D) -- (H_B);
\end{tikzpicture}
    \caption{Influence diagram for optimising the preventive care decision strategy for CHD.}
    \label{fig:CHD_influence_diagram}
\end{figure}

The problem setting is such that the patient is assumed to have a prior risk estimate $R_0$. A risk estimate is a prediction of the patient’s chance of having a CHD event in the next ten years. The risk estimates are grouped into risk levels, which range from 0\% to 100\% with a suitable discretisation, e.g., $S_{R_0}=\{0\%, 1\%, ..., 99\%, 100 \% \}$. We note that it might be beneficial to consider a less trivial discretisation that is finer in the region where most of the probability mass is assumed to lie and coarser elsewhere. Nevertheless, we chose to proceed as such since it requires no information on the probability distributions. The first testing decision $T_1$ is made based on the prior risk estimate. This entails deciding whether to perform TRS or GRS or if no testing is needed. If a test is conducted, the risk estimate is updated ($R_1$) and based on the new information a second testing decision $T_2$ follows. It entails deciding whether further testing should be conducted or not. The second testing decision is constrained so that the same test which was conducted in the first stage cannot be repeated. If a second test is conducted, the risk estimate is updated again ($R_2$). The treatment decision $T_D$ (dictating whether the patient receives preventive statin medicine or not) is made based on the resulting risk estimate of this testing process. Note that if no tests are conducted, the treatment decision is made based on the prior risk estimate. Figure \ref{fig:CHD_influence_diagram} provides an influence diagram for the decision problem. 

Node $H$ represents the uncertainty of whether the patient has a CHD event or remains healthy during the 10-year time frame. Node $H$ has the prior risk level $R_0$ in its information set because a premise of the modelling proposed in \citet{hynninen2019value} is that the prior risk accurately describes the
probability of having a CHD event, i.e., 
\begin{equation*}
    P(H = \text{CHD} \mid R_0=\alpha)=\alpha.
\end{equation*}
On the other hand, nodes $R_1$ and $R_2$ represent the updated risk level after the first and second test decisions, respectively. If a test is conducted, the risk estimate is updated using the Bayes' rule
\begin{equation*}
    P(\text{CHD} \mid \text{test result}) = \frac{ P(\text{test result} \mid \text{CHD}) \times P(\text{CHD}) }{P(\text{test result})},
\end{equation*}
where the conditional probabilities $P(\text{test result} \mid \text{CHD})$ are from \citet{abraham2016genomic} and the probability of having a CHD event, denoted by $P(\text{CHD})$, is the prior risk level $R_0$ or the updated risk level $R_1$, depending on whether it is the first or second test in question. The denominator $P(\text{test result})$ is calculated as a sum of the numerator and $P(\text{test result} \mid \text{no CHD}) \times  P(\text{no CHD})$, where $P(\text{no CHD}) = 1 - P(\text{CHD})$. As the states of nodes $R_i$, $i \in \{0,1,2\}$, represent risk levels, the probability of a state in these nodes is the probability of the given test updating the risk estimate to that level from the previous estimate.

The first and second testing decisions are represented by $T_1$ and $T_2$, respectively. Since conducting the same test twice is forbidden, all paths where the same test is repeated in $T_1$ and $T_2$ are included in the set of forbidden paths (cf. Section \ref{sec:formulations}). Furthermore, the forbidden paths include all paths where the first testing decision $T_1$ is to not perform testing but then the second testing decision $T_2$ is to perform a test. This is because the information yielded from performing only one test is not affected by whether the test is performed in the first or second stage of testing. Therefore, forbidding the paths where no test is performed in $T_1$ and a test is performed in $T_2$ reduces redundancy in the model without information loss. The final treatment decision is represented by node $T_D$, where the options are to provide or withhold treatment. The treatment decision is made based on the updated risk estimate represented by node $R_2$.

Since the first node in the influence diagram presented in Figure \ref{fig:CHD_influence_diagram} is the chance node $R_0$, any decision strategy would be conditioned on its realisation. This leads to a natural separability of the problem, meaning that it can be solved for individual risk levels $0\%, 1\%, \dots 100\%$. This has the benefit of allowing the calculations to be parallelised, at the expense of potentially causing inconsistencies related to e.g., multiple solutions in the MIP problem or rounding-induced errors.

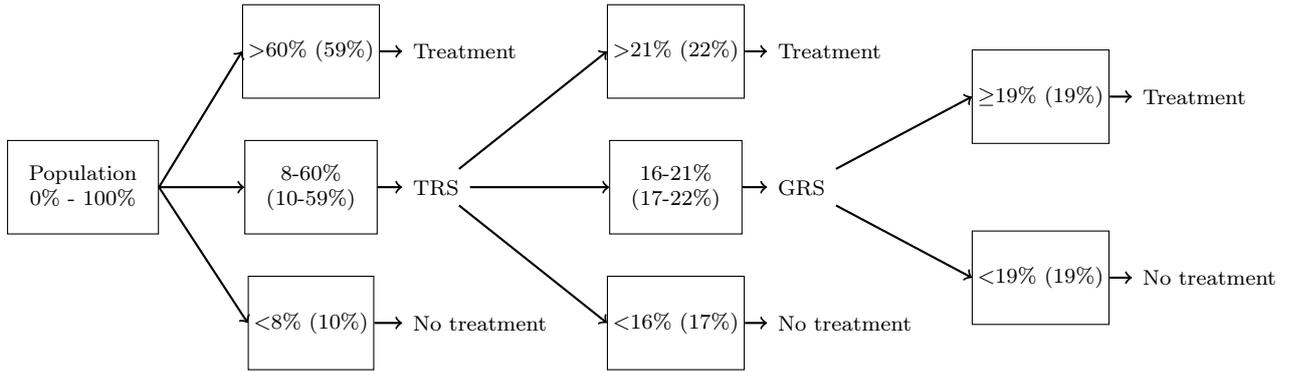
\begin{figure}
    \centering
    \begin{tikzpicture}
    [decision/.style={draw, minimum size=3em, inner sep=2pt}, 
    scale=1.2]
    \node[decision] (1) at (0, 2)  {\scriptsize \begin{tabular}{c} Population \\ 0\% - 100\% \end{tabular}};
    \node[decision] (2) at (2.5, 3.5)  {\scriptsize $>$60\% (59\%)};
    \node[decision] (3) at (2.5, 2)    {\scriptsize \begin{tabular}{c}8-60\% \\ (10-59\%)\end{tabular}};
    \node[decision] (4) at (2.5, 0.5)  {\scriptsize $<$8\% (10\%)};
    \node[right] (T1) at (3.5, 3.5) {\scriptsize Treatment};
    \node[right] (TRS) at (3.5, 2) {\scriptsize TRS};
    \node[right] (NT1) at (3.5, 0.5) {\scriptsize No treatment};
    \node[decision] (5) at (6.5, 3.5)  {\scriptsize $>$21\% (22\%)};
    \node[decision] (6) at (6.5, 2)    {\scriptsize \begin{tabular}{c}16-21\% \\ (17-22\%)\end{tabular}};
    \node[decision] (7) at (6.5, 0.5)  {\scriptsize $<$16\% (17\%)};
    \node[right] (T2) at (7.5, 3.5) {\scriptsize Treatment};
    \node[right] (GRS) at (7.5, 2) {\scriptsize GRS};
    \node[right] (NT2) at (7.5, 0.5) {\scriptsize No treatment};
    \node[decision] (8) at (10.5, 3)  {\scriptsize $\ge$19\% (19\%)};
    \node[decision] (9) at (10.5, 1)  {\scriptsize $<$19\% (19\%)};
    \node[right] (T3) at (11.5, 3) {\scriptsize Treatment};
    \node[right] (NT3) at (11.5, 1) {\scriptsize No treatment};
    \draw[thick, ->] (1.east) -- (2.west);
    \draw[thick, ->] (1.east) -- (3);
    \draw[thick, ->] (1.east) -- (4.west);
    \draw[thick, ->] (2) -- (T1);
    \draw[thick, ->] (3) -- (TRS);
    \draw[thick, ->] (4) -- (NT1);
    \draw[thick, ->] (TRS) -- (5.west);
    \draw[thick, ->] (TRS) -- (6);
    \draw[thick, ->] (TRS) -- (7.west);
    \draw[thick, ->] (5) -- (T2);
    \draw[thick, ->] (6) -- (GRS);
    \draw[thick, ->] (7) -- (NT2);
    \draw[thick, ->] (GRS) -- (8.west);
    \draw[thick, ->] (GRS) -- (9.west);
    \draw[thick, ->] (8) -- (T3);
    \draw[thick, ->] (9) -- (NT3);
\end{tikzpicture}
    \caption{Optimal strategy obtained by our model (in parentheses, the original value from \citet{hynninen2019value})}
    \label{fig:CHD_result}
\end{figure}

An interesting result is that the optimal strategy found by our model is the same strategy that was deemed the best among strategies (i-vi) in \citet{hynninen2019value}. In a way, this provides optimality guarantees to their results, which, in principle, they could not have determined without exhaustively testing all possible testing strategies. In addition, the optimal thresholds from our model correspond closely to those in \citet{hynninen2019value}. Figure \ref{fig:CHD_result} illustrates the strategy obtained by our model, indicating also the thresholds found in \citet{hynninen2019value} for comparison. We are confident that the small differences in the threshold values are simply artefacts related to the way the discretisation (i.e., rounding) was performed.

\begin{figure}[!htb]
     \centering
     \includegraphics[width=0.9\textwidth]{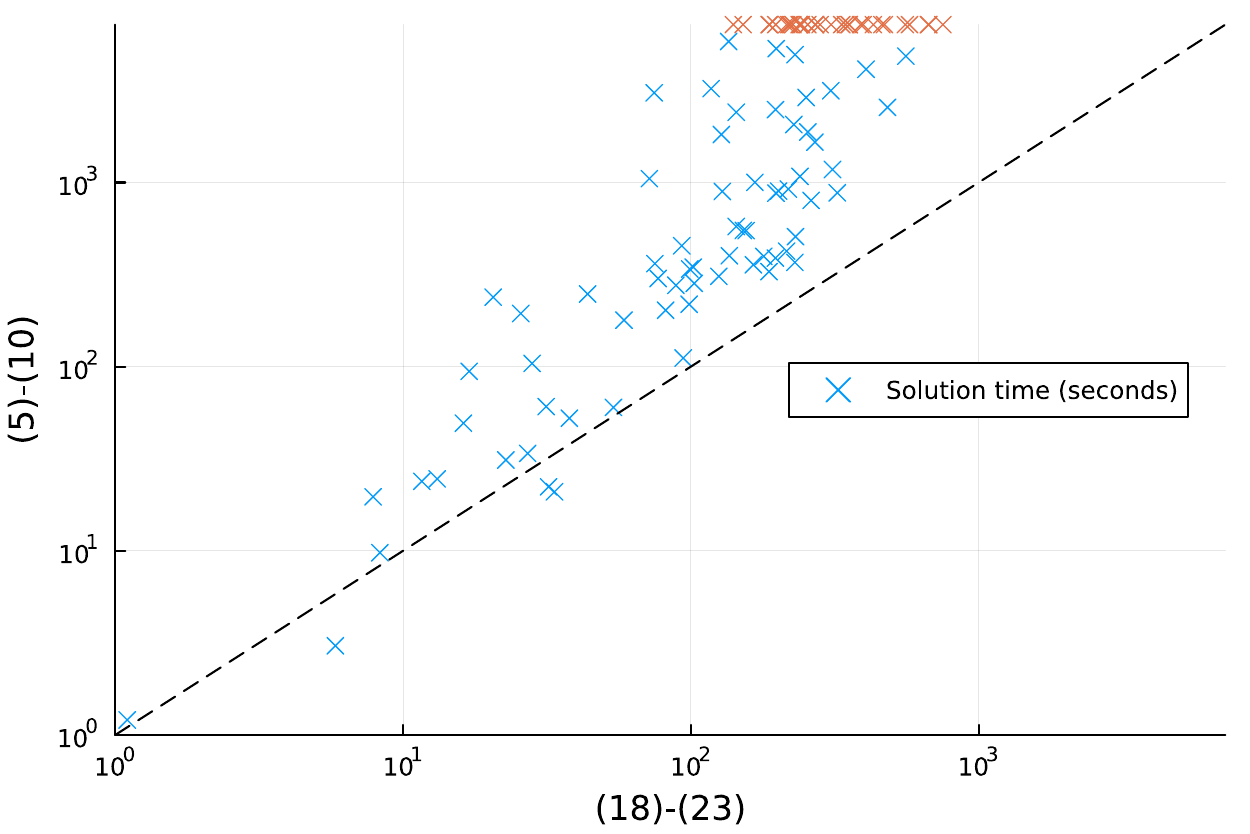}
     \caption{The solution times for each subproblem (prior risk level) in the CHD prevention case study with a time limit of 7200s. Points above the dashed line correspond to subproblems where the new formulation is faster than the old formulation \citep{salo2022} and orange points correspond to subproblems where the time limit was reached before finding the optimal solution.}
     \label{fig:CHD}
 \end{figure}
 
Finally, we continue the computational comparison from Section \ref{sec:computational_experiments} by solving the 101 subproblems corresponding to the prior risk levels 0\%, 1\%, ..., 100\% using the two formulations \eqref{eq:dp2_obj}-\eqref{eq:dp2_z_bin} and \eqref{eq:dp_obj}-\eqref{eq:dp_z_bin}. The results presented in Figure \ref{fig:CHD} make it evident that the new formulation proposed in this paper is significantly more efficient from a computational standpoint. Indeed, the formulation \eqref{eq:dp_obj}-\eqref{eq:dp_z_bin} from \citet{salo2022} could not find the optimal solution within 2h for almost a third of the subproblems (31 of 101). On the other hand, nearly all subproblems (98 of 101) were solved in under 10 minutes using the proposed formulation \eqref{eq:dp2_obj}-\eqref{eq:dp2_z_bin}.
\section{Conclusions} \label{sec:conclusions}

In this paper, we expand on the ideas originally proposed in \citet{salo2022} providing multiple methodological enhancements. These enhancements include a novel and more efficient formulation, valid bounds to tighten relaxations, and a heuristic which can be used to find feasible solutions and, consequently, to warm start the MIP solver. 
Furthermore, we conduct a novel case study based on the study originally proposed by \citet{hynninen2019value}. Our objective is to demonstrate that the proposed models can be used in settings which would normally require resorting to more ad-hoc computational tools, lending themselves to be a general and accessible tool for practitioners. We believe that this will allow for a much wider range of practitioners and researchers to have access to mathematical optimisation-based tools for supporting decision-making. Furthermore, this will create novel inroads for the use of mathematical optimisation in the area of decision analysis at large, potentially unveiling new and promising directions for future developments.

In terms of alternative further developments, we see several directions that deserve further investigation. First, decision programming as a modelling framework is still in its infancy, and, consequently, many obstacles are still to be overcome for its widespread adoption. One of these obstacles is computational requirements. Decision programming models grow large as the number of nodes and/or states increase, and thus it would greatly benefit from alternative ideas that can tackle such large-scale problems. These can be, for example, related to alternative formulations that convert the influence diagram into an intermediate structure and employ ideas from Bayesian inference to yield a more compact MIP model (see \citet{parmentier2020integer}). Another direction worth exploring is the employment of decomposition methods, in particular, those which allow for a delayed generation of structural elements of the model, in our case the paths $s \in S$ (see Section \ref{sec:decision_programming}). Another interesting avenue would be to pursue methods that can reap benefits from employing parallelization, given the increasing availability of high-performance computing clusters.

\section*{Acknowledgements}

We are enormously grateful for the input from Juho Andelmin, whose initial implementations led to the development of \texttt{DecisionProgramming.jl}. We are also grateful for the contributions of a number of graduate and undergraduate students to the development of the package, as well as to the welcoming and supportive JuMP community. We also gratefully acknowledge the financial support from the Research Council of Finland (decision number 332180). Finally, the computational experiments were performed using computer resources within the Aalto University School of Science “Science-IT” project.

\bibliographystyle{plainnat}
\bibliography{references}

\begin{thebibliography}{16}
\providecommand{\natexlab}[1]{#1}
\providecommand{\url}[1]{\texttt{#1}}
\expandafter\ifx\csname urlstyle\endcsname\relax
  \providecommand{\doi}[1]{doi: #1}\else
  \providecommand{\doi}{doi: \begingroup \urlstyle{rm}\Url}\fi

\bibitem[Abraham et~al.(2016)Abraham, Havulinna, Bhalala, Byars, De~Livera,
  Yetukuri, Tikkanen, Perola, Schunkert, Sijbrands, et~al.]{abraham2016genomic}
Gad Abraham, Aki~S Havulinna, Oneil~G Bhalala, Sean~G Byars, Alysha~M
  De~Livera, Laxman Yetukuri, Emmi Tikkanen, Markus Perola, Heribert Schunkert,
  Eric~J Sijbrands, et~al.
\newblock Genomic prediction of coronary heart disease.
\newblock \emph{European Heart Journal}, 37\penalty0 (43):\penalty0 3267--3278,
  2016.

\bibitem[Bertsimas and Dunn(2019)]{bertsimas2019machine}
Dimitris Bertsimas and Jack Dunn.
\newblock \emph{Machine learning under a modern optimization lens}.
\newblock Dynamic Ideas LLC Charlestown, MA, 2019.

\bibitem[Bielza et~al.(2000)Bielza, G{\'o}mez, R{\i}os-Insua, and del
  Pozo]{bielza2000structural}
Concha Bielza, Manuel G{\'o}mez, S~R{\i}os-Insua, and JA~Fern{\'a}ndez del
  Pozo.
\newblock Structural, elicitation and computational issues faced when solving
  complex decision making problems with influence diagrams.
\newblock \emph{Computers \& Operations Research}, 27\penalty0 (7-8):\penalty0
  725--740, 2000.

\bibitem[Bielza et~al.(2011)Bielza, G{\'o}mez, and Shenoy]{bielza2011review}
Concha Bielza, Manuel G{\'o}mez, and Prakash~P Shenoy.
\newblock A review of representation issues and modeling challenges with
  influence diagrams.
\newblock \emph{Omega}, 39\penalty0 (3):\penalty0 227--241, 2011.

\bibitem[Birge and Louveaux(2011)]{birge2011introduction}
John~R Birge and Francois Louveaux.
\newblock \emph{Introduction to stochastic programming}.
\newblock Springer Science \& Business Media, 2011.

\bibitem[Cobb(2024)]{cobb2024intermittent}
Barry~R Cobb.
\newblock Intermittent sampling for statistical process control with the number
  of defectives.
\newblock \emph{Computers \& Operations Research}, 161:\penalty0 106423, 2024.

\bibitem[Howard and Matheson(2005)]{howard2005influence}
Ronald~A Howard and James~E Matheson.
\newblock Influence diagrams.
\newblock \emph{Decision Analysis}, 2\penalty0 (3):\penalty0 127--143, 2005.

\bibitem[Hynninen et~al.(2019)Hynninen, Linna, and
  Vilkkumaa]{hynninen2019value}
Yrj{\"a}n{\"a} Hynninen, Miika Linna, and Eeva Vilkkumaa.
\newblock Value of genetic testing in the prevention of coronary heart disease
  events.
\newblock \emph{PloS one}, 14\penalty0 (1):\penalty0 e0210010, 2019.

\bibitem[Lauritzen and Nilsson(2001)]{lauritzen2001representing}
Steffen~L Lauritzen and Dennis Nilsson.
\newblock Representing and solving decision problems with limited information.
\newblock \emph{Management Science}, 47\penalty0 (9):\penalty0 1235--1251,
  2001.

\bibitem[Mau{\'a} et~al.(2012)Mau{\'a}, de~Campos, and
  Zaffalon]{maua2012solving}
Denis~Deratani Mau{\'a}, Cassio~P de~Campos, and Marco Zaffalon.
\newblock Solving limited memory influence diagrams.
\newblock \emph{Journal of Artificial Intelligence Research}, 44:\penalty0
  97--140, 2012.

\bibitem[Parmentier et~al.(2020)Parmentier, Cohen, Lecl{\`e}re, Obozinski, and
  Salmon]{parmentier2020integer}
Axel Parmentier, Victor Cohen, Vincent Lecl{\`e}re, Guillaume Obozinski, and
  Joseph Salmon.
\newblock Integer programming on the junction tree polytope for influence
  diagrams.
\newblock \emph{INFORMS Journal on Optimization}, 2\penalty0 (3):\penalty0
  209--228, 2020.

\bibitem[Pauker and Kassirer(1980)]{pauker1980threshold}
Stephen~G Pauker and Jerome~P Kassirer.
\newblock The threshold approach to clinical decision making.
\newblock \emph{New England Journal of Medicine}, 302\penalty0 (20):\penalty0
  1109--1117, 1980.

\bibitem[Puterman(1990)]{puterman1990markov}
Martin~L Puterman.
\newblock Markov decision processes.
\newblock \emph{Handbooks in operations research and management science},
  2:\penalty0 331--434, 1990.

\bibitem[Salo et~al.(2022)Salo, Andelmin, and Oliveira]{salo2022}
Ahti Salo, Juho Andelmin, and Fabricio Oliveira.
\newblock Decision programming for mixed-integer multi-stage optimization under
  uncertainty.
\newblock \emph{European Journal of Operational Research}, 299\penalty0
  (2):\penalty0 550--565, 2022.

\bibitem[Schrijver et~al.(2003)]{schrijver2003combinatorial}
Alexander Schrijver et~al.
\newblock \emph{Combinatorial optimization: polyhedra and efficiency},
  volume~A.
\newblock Springer, 2003.

\bibitem[Tikkanen et~al.(2013)Tikkanen, Havulinna, Palotie, Salomaa, and
  Ripatti]{tikkanen2013genetic}
Emmi Tikkanen, Aki~S Havulinna, Aarno Palotie, Veikko Salomaa, and Samuli
  Ripatti.
\newblock Genetic risk prediction and a 2-stage risk screening strategy for
  coronary heart disease.
\newblock \emph{Arteriosclerosis, thrombosis, and vascular biology},
  33\penalty0 (9):\penalty0 2261--2266, 2013.

\end{thebibliography}

\end{document}